\input amstex\documentstyle{amsppt}  
\pagewidth{12.5cm}\pageheight{19cm}\magnification\magstep1
\topmatter
\title {Algebraic and geometric methods in representation theory}\endtitle
\author G. Lusztig\endauthor
\thanks{Supported in part by National Science Foundation grant DMS-1303060 and by a Simons Fellowship.
This paper will appear as part of a book to be published by the Shaw Prize Foundation.}
\endthanks
\address{Department of Mathematics, M.I.T., Cambridge, MA 02139}\endaddress
\endtopmatter   
\document

\define\Irr{\text{\rm Irr}}

\define\ul{\un l}

\define\bo{\bar o}

\define\bY{\bar Y}

\define\bco{\bar{\co}}

\define\sqc{\sqcup}

\define\bX{\bar X}

\define\op{\oplus}
   
\define\part{\partial}
\define\emp{\emptyset}

\define\m{\mapsto}
\define\do{\dots}

\define\bsl{\backslash}

\define\sub{\subset}    

\define\T{\times}
\define\ti{\tilde}
\define\nl{\newline}
\redefine\i{^{-1}}
\define\fra{\frac}
\define\un{\underline}
\define\ov{\overline}
\define\ot{\otimes}
\define\bbq{\bar{\QQ}_l}

\define\tr{\text{\rm tr}}

\redefine\c{\chi}
\define\g{\gamma}

\define\io{\iota}

\define\p{\pi}
\define\ph{\phi}

\define\r{\rho}
\define\s{\sigma}
\redefine\t{\tau}
\define\th{\theta}

\redefine\l{\lambda}
\define\z{\zeta}

\redefine\G{\Gamma}

\define\Ph{\Phi}

\define\CC{\bold C}

\define\FF{\bold F}
\define\GG{\bold G}

\define\NN{\bold N}

\define\QQ{\bold Q}

\define\ZZ{\bold Z}

\define\ca{\Cal A}
\define\cb{\Cal B}

\define\cf{\Cal F}

\define\ck{\Cal K}

\define\co{\Cal O}

\define\ct{\Cal T}

\define\fS{\frak S}

\define\fU{\frak U}

\define\tW{\ti W}
\define\tX{\ti X}

\subhead 1\endsubhead
Being awarded the Shaw Prize was a (pleasant) surprise for me; I feel very honored by it.

This paper is an expanded version of the Shaw Prize Lecture given on September 25, 2014 at the Chinese 
University of Hong Kong; it is concerned with some of my work in the theory of group representations, with 
emphasis on representations of finite groups of Lie type. Throughout this paper I have tried to explain not 
only what the results are but also the origin of those results.

I wish to thank David Vogan for his comments on an earlier version of this paper.

\subhead 2\endsubhead
In mathematics, groups are everywhere. They also appear in the natural sciences: in physics, the motions of
space-time form a (continuous) group; in chemistry, the symmetries of a crystal form a (finite) group. We
shall be mainly interested in finite groups and more specifically, in finite simple groups, which form the 
building blocks for all finite groups. Rather surprisingly, up to finitely many known exceptions, and up to 
extensions by finite commutative groups, the finite noncommutative simple groups are of Lie type that is,
they can be described in terms of continuous (or Lie) groups as follows.

Start with a compact connected Lie group; by complexification (Chevalley, 1946) this gives rise to a complex
(reductive connected) Lie group and by a careful choice of integral structure (Chevalley, 1955 and 1960), to
a reductive connected algebraic group $G$ over any given field, for example over $K$, an algebraic closure 
of the finite field with $p$ elements (with $p$ a prime number). For any field $k$ we shall denote by 
$GL_n(k)$ the group of automorphisms of $k^n$ (or the group of invertible $n\T n$ matrices with entries in 
$k$); let $F_1:GL_n(K)@>>>GL_n(K)$ be the homomorphism given by $(a_{ij})\m(a_{ij}^p)$. We now choose an 
endomorphism $F:G@>>>G$ such that for some $n,s,s'$ in $\ZZ_{>0}$ and some imbedding $\io:G@>>>GL_n(K)$ as a
closed subgroup, we have $\io(F^{s'}(g))=F_1^s(\io(g))$ for any $g\in G$. Then the fixed point set $G^F$ of 
$F$ is a finite group, said to be of Lie type. When $s'=1$, $G^F$ is the group $G(\FF_q)$ for an 
$\FF_q$-rational structure on $G$ with Frobenius map $F$ (here $\FF_q$ is the subfield of $K$ with $q=p^s$ 
elements). The Weyl group $W$ of $G$ should be also regarded as a finite group of Lie type. 

For example, when $G=GL_n(K)$, the finite group $GL_n(\FF_q)$ is of Lie type; $W$ can be identified with the 
group of all permutations of $\{1,2,\do,n\}$. As another example, let $G=Sp_{2n}(K)$ be the group of all 
$g\in GL_{2n}(K)$ which preserve the symplectic form 
$$x_1y_{2n}-x_{2n}y_1+x_2y_{2n-1}-x_{2n-1}y_2+\do+x_ny_{n+1}-x_{n+1}y_n$$ 
(a symplectic group) and let $F:G@>>>G$ be the restriction of $F_1^s:GL_{2n}(K)@>>>F_1:GL_{2n}(K)$. Then
$G^F=Sp_{2n}(\FF_q)$ is a finite group of Lie type (here $q=p^s$). In this case $W$ can be identified with 
the group of all permutations of $\{1,2,\do,2n\}$ which commute with the involution 
$1\m2n,2\m2n-1,3\m3n-2,\do,2n\m1$.

\subhead 3\endsubhead
A representation of a finite group $\G$ is a homomorphism $\r:\G@>>>GL(V)$ of $\G$ into the group of 
automorphisms of a $\CC$-vector space $V$ (all vector spaces in this paper are assumed to be of finite 
dimension). We say that $(\r,V)$ is {\it irreducible} if $V\ne0$ and there is no vector subspace $V'$ of $V$
(with $0\ne V'\ne V$) such that $\r(g)V'=V'$ for any $g\in\G$. (We will occasionally replace $\CC$ by other 
fields.) These definitions appeared in a groundbreaking paper of F. G. Frobenius in 1896 which marks the 
birth of representation theory. The study of irreducible representations of a finite group is a key to 
understanding the finite group itself in the same way as understanding an object can be achieved by 
analyzing pictures of that object from many different angles. 

For any representation $(\r,V)$ of $\G$, the character $\c_\r$ of $\r$ is the function $\G@>>>\CC$ given by 
$g\m\tr(\r(g),V)$. We have $\c_\r(g)=\c_r(g'gg'{}\i)$ for any $g,g'\in\G$; thus $\c_\r$ is constant on each
conjugacy class of $\G$. 

Let $\Irr\G$ be the set of isomorphism classes of irreducible representations of $\G$. In his 1896 paper, 
Frobenius showed that, if $\Irr\G=\{\r_1,\r_2,\do,\r_e\}$, then $e$ is also the number of conjugacy classes 
of $\G$ and the functions $\c_{\r_1},\c_{\r_2},\do,\c_{\r_e}$ (the ``irreducible characters'') form a basis 
for the $\CC$-vector space of functions $\G@>>>\CC$ which are constant on the conjugacy classes 
$C_1,C_2,\do,C_e$ of $\G$. (If $g_i\in C_i$, then the invertible $e\T e$-matrix $(\c_{\r_j}(g_i))$ is called
the {\it character table} of $\G$.) In the same paper Frobenius defined the notion of representation of $\G$
induced by a representation of a subgroup and several years later, in 1900, he showed that the irreducible 
characters of the symmetric group in $n$ letters can be expressed as explicit $\ZZ$-linear combinations of 
characters of representations induced by the unit representations of various subgroups.

In this paper we are interested in understanding as much as possible about the representations of $G^F$ where
$G$, a connected reductive group over $K$, and $F:G@>>>G$ (as in no.2) are fixed; we assume that $F$ is the 
Frobenius map for an $\FF_q$-rational structure on $G$ where $q=p^s$ for some $s\in\ZZ_{>0}$ so that 
$G^F=G(\FF_q)$. As in no.2, $W$ will denote the Weyl group of $G$.

\subhead 4\endsubhead
In 1955, J. A. Green published a remarkable paper in which he described completely the character table of 
the general linear groups $GL_n(\FF_q)$. (For $n=2$ this was essentially done by Frobenius in 1896.) Green 
used Frobenius's method of taking induced representations from proper subgroups (in this case the subgroups 
are stabilizers of flags of subspaces in $\FF_q^n$) and some rather nontrivial combinatorics coming from 
Hall algebras. But then he was faced with the problem of constructing the ``discrete series'', a family of 
irreducible representations which are not contained in any of the induced representations above and for 
which a new method was needed. His approach was based on the following result that he proved using 
R. Brauer's characterization of characters.

(a) {\it Let $u$ be an isomorphism of the group $K^*$ with the group of roots of $1$ of order prime to $q$ 
in $\CC$. Let $\G$ be a finite group and let $\ph:\G@>>>GL_N(K)$ be a group homomorphism. Define a function 
$\c:\G@>>>\CC$ by $\c(g)=u(\l_1)+\do+u(\l_N)$ where $\l_1,\do,\l_N$ are the eigenvalues (in $K^*$) of 
$\ph(g):K^N@>>>K^N$; here $g\in\G$. Then $\c$ is a $\ZZ$-linear combination of irreducible characters of 
$\G$.}
\nl
(Green called $\c$ the ``Brauer lifting'' of $\ph$.) (a) allowed Green to find the characters of the 
discrete series representations but not the representations themselves. It turns out that they are all of 
dimension $(q-1)(q^2-1)\do(q^{n-1}-1)$.

It is interesting that Green's work on the representations of $GL_n(\FF_q)$ was almost at the same time as 
the work of Harish-Chandra on representations of real semisimple Lie groups. In both cases, most irreducible
representations were associated with characters of maximal tori defined over the ground field and in both 
cases there was at most one ``elliptic'' maximal torus up to conjugacy.

The next progress in this field was achieved in 1968, when B.Srinivasan computed the character table of the 
symplectic group $Sp_4(\FF_q)$ for $q$ odd. In this case there are two conjugacy clases of elliptic maximal 
tori and an unexpected irreducible discrete series representation appears.

\subhead 5\endsubhead
In 1970, at the Institute for Advanced Study, I attended a talk by D. Quillen, in which he explained his 
solution to a problem in homotopy theory, the Adams conjecture. That solution used essentially the Brauer 
lifting 4(a) of the natural $n$-dimensional representation of $GL_n(\FF_q)$. After the talk, I have asked 
M. Atiyah whether the irreducible representations which enter in that Brauer lifting were explicitly known, 
and he told me that only their characters were known. I became very interested to understand the 
representations which enter, not only their characters. In 1971 I accepted a lectureship at the University 
of Warwick where J. A. Green was a professor. At Warwick, I learned much about algebraic groups from 
discussions with R. W. Carter with whom I also wrote two papers. In one of those papers I became familiar 
with Hecke algebras (with parameter $0$) and how to use them to construct irreducible modular 
representations. The experience with Hecke algebras was to become very useful in my later research. 

In 1972 and early 1973, I found a way to describe explicitly the Brauer lifting of the natural 
$n$-dimensional representation of $GL_n(\FF_q)$. I will try to explain some of the ideas involved, assuming 
for simplicity that $q=p$ is a prime number and that $\CC$ is replaced by a maximal unramified extension of 
$\QQ_p$, the $p$-adic numbers (in this case there is a natural choice for $u$ in 4(a) whose image is 
contained in the units of $\ZZ_p$, the $p$-adic integers).

Let $V$ be an $n$-dimensional $\FF_p$-vector space, $n\ge1$. Let $E$ be the set of complete flags
$V_0\sub V_1\sub\do\sub V_n$ in $V$ ($\dim V_i=i$ for all $i$). Let $X_V$ be the $\FF_p$-vector space 
consisting all functions which to each $(V_0\sub V_1\sub\do\sub V_n)\in E$ associate a vector 
$f(V_0\sub V_1\sub\do\sub V_n)\in V_1$ in such a way that for any $i\in[1,n-1]$ and any almost complete flag
of the form $(V_0\sub V_1\sub\do\sub V_{i-1}\sub V_{i+1}\sub\do\sub V_n)$, the sum 
$$\sum_{V_i}f(V_0\sub V_1\sub\do\sub V_{i-1}\sub V_i\sub V_{i+1}\sub\do\sub V_n)$$
is $0$ (here $V_i$ runs over all $i$-dimensional subspaces which fit between $V_{i-1}$ and $V_{i+1}$); the 
sum is taken in $V_1$ if $i\ge2$ and in $V_2$ if $i=1$. Note that $X_V$ has a natural action of $GL(V)$. (If
in the definition of $X_V$ one takes functions $f$ with values in $\CC$ instead of $V_1$ and satisfying the 
similar equations, one obtains the Steinberg representation of $GL(V)$ of dimension $p^{1+2+\do+(n-1)}$.)
One shows that $\dim X_V=(p-1)(p^2-1)\do(p^{n-1}-1)$ and that $X_V$ is naturally the reduction modulo $p$ of
a free $\ZZ_p$-module $\tX_V$ of rank $(p-1)(p^2-1)\do(p^{n-1}-1)$ with natural action of $GL(V)$ (this is 
the most nontrivial part of the story). Moreover, if for any $j\in[1,n]$ we set $\tX_V^j=\op_{V_j}\tX_{V_j}$
(sum over the $j$-dimensional subspaces of $V$) then $\tX_V^j\ot\QQ_p$ has a natural action of $GL(V)$ and 
these provide the required components of the Brauer lifting, which is equal to 
$\sum_{j=1}^n(-1)^{j-1}\c_{\tX_V^j}$.

In late 1973 I found the explicit structure of the Brauer lifting of the standard representation of the
various classical groups over $\FF_q$, and I thus obtained in each case some discrete series representations
which were new (at that time). But it seemed to be difficult to obtain other discrete series by this method. 

\subhead 6\endsubhead
In the early 1960's A. Grothendieck defined $l$-adic cohomology spaces for algebraic varieties in 
characteristic $p$ (here $l$ is a prime number other than $p$) which shared many properties with the usual 
rational cohomology of algebraic varieties over $\CC$. For example, the smooth curve $C_0$ with equation  
$x^{q+1}+y^{q+1}+z^{q+1}=0$ in the projective plane over $K$ has $l$-adic cohomology group in degree $1$ of 
dimension $q^2-q$, the same as it would have (in ordinary cohomology) if it was viewed as a curve over 
$\CC$. But the curve over $\CC$ has no interesting automorphism, while $C_0$ has a natural action of a 
finite unitary group (a subgroup of $GL_3(\FF_{q^2})$) which, as Tate and Thompson have shown in a 1965 
paper (by J. Tate), induces an irreducible representation of dimension $q^2-q$ of this unitary group on the 
first $l$-adic cohomology space. (This representation is difficult to obtain by other means.)

Towards the end of 1973, T. A. Springer explained to me some work of V. Drinfeld about which he learned 
during a recent visit to Moscow. Namely, Drinfeld (who was 19 years old at the time) considered the plane 
curve $C=\{(x,y)\in K^2;x^qy-xy^q=1\}$ on which $SL_2(\FF_q)$ acts naturally (an old observation of 
L. E. Dickson) and $\ct:=\{\l\in K;\l^{q+1}=1\}$ acts (freely) by homothety, commuting with the 
$SL_2(\FF_q)$-action. 
(Note that $SL_2(\FF_q)=G^F$ in the case where $G=SL_2(K)$ and $F:G@>>>G$ is given by 
$F(a_{ij})=(a_{ij}^q)$.) Then $SL_2(\FF_q)\T\ct$ acts on $H^1_c(C,\bbq)$, the $l$-adic cohomology with 
compact support in degree $1$ of $C$. (Here $\bbq$ is an algebraic closure of the field of $l$-adic 
numbers.) Note that $SL_2(\FF_q)\T\ct$ can be viewed as a subgroup of the unitary group above and $C$ can be
viewed as the complement of a finite set in $C_0$ above so that the action of $SL_2(\FF_q)\T\ct$ on $C$ 
becomes the restriction of the action of the unitary group on $C_0$. Drinfeld observed that all irreducible 
$(q-1)$-dimensional representations of $SL_2(\FF_q)$ appear in the action on $H^1_c(C,\bbq)$ (as eigenspaces
of the $\ct$-action); these are restrictions of discrete series representations of $GL_2(\FF_q)$. (Most of 
the remaining irreducible representation have dimension $q+1$ and are easy to obtain by induction from a 
subgroup.) 

\subhead 7\endsubhead
During the spring of 1974 I was at the Institut des Hautes \'Etudes Scientifiques in Bures-sur-Yvette and 
in the later part of my stay I started a joint work with P. Deligne, trying to generalize Drinfeld's 
construction in no.6 to a general $G^F$; this was 
written up as a joint paper in the first half of 1975 (it appeared in 1976). 

In the setup of no.6, the map $(x,y)\m[\text{ line in $K^2$ spanned by }(x,y)]$ identifies the quotient 
$\ct\bsl C$ with $\{L\in\cb;F(L)\ne L\}$ where $\cb$ denotes the set of lines in $K^2$ and $F:\cb@>>>\cb$ is
induced by $(x,y)\m(x^q,y^q)$. Thus $\cb$ has a natural partition 
$\{L\in\cb;F(L)=L\}\sqc\{L\in\cb;F(L)\ne L\}$ in which the second piece is responsible for the
irreducible representations of dimension $q-1$ and similarly the first piece is responsible for the
irreducible representations of dimension $q+1$. (Note also that the first piece has Euler characteristic
$q+1$ and the second piece has Euler characteristic $1-q$.)

In the general case, $\cb$ can be interpreted as the variety of all Borel subgroups (maximal closed 
connected solvable subgroups) of $G$ (when $G=SL_2(K)$, a line $L$ in $K^2$ can be identified with its 
stabilizer of $L$ in $SL_2(K)$). Now $G$ acts on $\cb\T\cb$ by simultaneous conjugation; it is known that
this action has finitely many orbits naturally indexed by the elements of $W$ (when $G=SL_2(K)$ there are 
two such orbits: one, indexed by the unit element, is the diagonal and the other, indexed by the element of 
order $2$, is the complement of the diagonal). See also no.8. We write $\co_w$ for the orbit indexed by 
$w\in W$. Now $F:G@>>>G$ induces an endomorphism of $\cb$ denoted again by $F$.
For any $w\in W$ we set $X_w=\{B\in\cb;(B,F(B))\in\co_w\}$, a subvariety of $\cb$. The subvarieties $X_w$ 
($w\in W$) form a partition of $\cb$ which in the case where $G=SL_2(K)$ coincides with the partition of 
$\cb$ described earlier. Each $X_w$ is stable under the conjugation action of $G(\FF_q)$. 
Now for each $w\in W$ we can find an $F$-stable maximal torus $\ct_w$ of $G$ which is contained in some
$B\in X_w$. Moreover, $X_w$ has a natural finite principal covering 
$\tX_w@>>>X_w$ with group $\ct_w^F=\ct_w\cap G^F$, to which the $G^F$-action can be lifted. (In the case 
where $G=SL_2(K)$ and $w\ne1$, $\tX_w$ can be taken to be the curve $C$ considered by Drinfeld.) It follows 
that $G^F$ acts naturally on the $l$-adic cohomology spaces $H^i_c(\tX_w,\bbq)$. The (commutative) group 
$\ct_w^F$ also acts naturally on $H^i_c(\tX_w,\bbq)$, commuting with the $G^F$-action; hence for any 
homomorphism $\th:\ct_w^F@>>>\bbq^*$, $G^F$ acts naturally on the $\th$-eigenspace $H^i_c(\tX_w,\bbq)_\th$ 
of the $\ct_w^F$-action. In the paper with Deligne we proved that any irreducible representation $\r$ of 
$G^F$ appears in the virtual representation $\sum_i(-1)^iH^i_c(\tX_w,\bbq)_\th$ for some $w,\th$. Now from 
each pair $(w,\th)$ as above one can produce in a natural way a conjugacy class $\g_{w,\th}$ defined over 
$\FF_q$ of elements of order prime to $q$ in a connected reductive group $G^*$ over $K$ (of Langlands dual 
type to that of $G$) with a natural $\FF_q$-structure and in the paper we prove that if $w,\th$ corresponds 
to $\r$ as above, then $\g_{w,\th}$ depends only on $\r$, not on $w,\th$ hence it can be denoted by $\g_\r$.
Thus we obtain a natural (surjective) map $\Ph:\r\m\g_\r$ from the set of irreducible representations of 
$G^F$ (up to isomorphism) to the set of $G^*$-conjugacy classes (defined over $\FF_q$) of elements of order 
prime to $q$ in $G^*$. The irreducible representations of $G^F$ in $\Ph\i(1)$ are called the unipotent 
representations. They are precisely the irreducible representations which appear in $H^i_c(X_w,\bbq)$ for 
some $i,w$ or, equivalently, which appear in the virtual representation $\sum_i(-1)^iH^i_c(X_w,\bbq)$ for 
some $w$. At the other extreme, for almost any $G^*$-conjugacy class $\g$ defined over $\FF_q$ of elements 
of order prime to $q$ in $G^*$, $\Ph\i(\g)$ is a single irreducible representation of $G^F$; it is equal to
$\pm\sum_i(-1)^iH^i_c(\tX_w,\bbq)_\th$ for some $w$ and some $\th$. Almost all irreducible representations 
of $G^F$ are obtained in this way. (In particular this proved a conjecture formulated by I. G. Macdonald in 
1968.)

\subhead 8\endsubhead
I will now describe in more detail how the $G$-orbits on $\cb\T\cb$ are parametrized by $W$ when $G$ is 
$GL_n(K)$ or $Sp_{2n}(K)$.

Assume first that $G=GL_n(K)$; then $\cb$ can be identified with the set of complete flags of subspaces 
$V_0\sub V_1\sub\do\sub V_n$ in $K^n$ ($\dim V_i=i$ for all $i$). Given two complete flags 
$V_0\sub V_1\sub\do\sub V_n$ and $V'_0\sub V'_1\sub\do\sub V'_n$, we can define uniquely a permutation 
$w$ of $\{1,2,\do,n\}$ by the following requirement: there exists a basis
$\{v_1,v_2,\do,v_n\}$ of $K^n$ such that for $i=1,2,\do,n$, $\{v_1,v_2,\do,v_i\}$ is a basis of $V_i$ and
$\{v_{w(1)},v_{w(2)},\do,v_{w(i)}\}$ is a basis of $V'_i$. The set of pairs of complete flags whose
associated permutation is $w$ form the orbit $\co_w$.
 
Assume next that $G=Sp_{2n}(K)$; then $\cb$ can be identified with the set of complete flags of subspaces 
$V_0\sub V_1\sub\do\sub V_{2n}$ in $K^{2n}$ ($\dim V_i=i$ for all $i$) such that for any $i=0,1,\do,2n$, 
$V_i$ is the perpendicular to $V_{2n-i}$ with respect to the symplectic form. To two such complete flags 
$V_0\sub V_1\sub\do\sub V_{2n}$ and $V'_0\sub V'_1\sub\do\sub V'_{2n}$, we associate a permutation $w$ of 
$\{1,2,\do,2n\}$ as for $GL_{2n}(K)$; this permutation commutes with the involution 
$1\m2n,2\m2n-1,3\m3n-2,\do,2n\m1$ hence is an element of $W$. The set of pairs of complete flags as above 
whose associated permutation is $w$ form the orbit $\co_w$.

\subhead 9\endsubhead
In late 1975, I wrote a paper (which appeared in 1976) in which I analyzed the variety $X_w$ in no.7 in the 
case where $G$ modulo its centre is simple and $X_w$ is irreducible, of minimum possible dimension. For 
simplicity we shall assume here that $G$ is split over $\FF_q$ that is, there exists a maximal torus $T$ of 
$G$ such that $F(t)=t^q$ for all $t\in T$. In this case, $w$ is a ``Coxeter element of minimal length'' of 
$W$. For example when $G=GL_n(K)$, one can take $X_w$ to be the set of complete flags 
$V_0\sub V_1\sub\do\sub V_n$ in $K^n$ (see no.8) such that $V_i\ne F(V_i)\sub V_{i+1}$ for $i=1,2,\do,n-1$. 
(Here $F:K^n@>>>K^n$ is given by $(x_i)\m(x_i^q)$.) In the case of
classical groups, the corresponding series of representations is the discrete series which can be approached
using Brauer lifting (see the last paragraph in no.5).
Note that in our case, $X_y$ is stable under $F:\cb@>>>\cb$ for any $y\in W$.
The main result of the paper is that the Frobenius map acts on $\op_iH^i_c(X_w,\bbq)$ semisimply and that
its eigenspaces are distinct irreducible representations of $G^F$. When $G$ is of exceptional type
one finds several (unipotent) discrete series representations among these eigenspaces. Moreover, the 
eigenvalues of the Frobenius map are explicitly computed. They are most of the time roots of $1$ times 
integer powers of $q$, but in some cases (type $E_7$, $E_8$) they can be odd powers of $\sqrt{-q}$. This 
shows that the topology of $X_w$ can be quite complicated and that I was quite lucky that these eigenvalues 
could be computed at all.
The results in this paper (in the case of exceptional groups) played an essential role in my later 
classification of unipotent representations of $G^F$.

\subhead 10\endsubhead
In the summer and fall of 1976 I wrote a paper (published in 1977) in which I found a complete 
classification (including dimensions) of the irreducible representations of $G^F$ in the case where $G$ is 
a classical group with connected centre, other than $GL_n(K)$. (Thus, for example, instead of $Sp_{2n}(K)$, 
I considered the group of symplectic similitudes.) 

The method I used was based on an extension of the method in no.7, which allowed cohomological induction 
from subgroups more general than a maximal torus, and on the use of the dimension formulas for the 
irreducible representations of Hecke algebras of type B with two parameters due to Hoefsmit (1974). This 
paper establishes what in no.17 is called ``quasi-induction'' for the representations of classical groups 
(with connected centre). It also establishes the parametrization of unipotent representations for these 
groups in terms of some new combinatorial objects, the ``symbols'', and the
classification of unipotent discrete series representations of classical
groups. Moreover, it is shown that the endomorphism algebra of the representation 
induced from an isolated discrete series representation 
to a larger classical group is an Iwahori-Hecke algebra (anticipating a later result of 
Howlett-Lehrer, 1980) and giving also precise information on the values of 
the parameters of that Iwahori-Hecke algebra. For this we need to count in terms of 
generating functions the number of conjugacy classes in a classical group with connected 
centre. This together with an inductive hypothesis and the methods outlined above give a
way to predict the number of isolated discrete series representations. The degrees of these 
representations can be guessed using the technique of symbols by ``interpolation'' from 
the degrees of already known representations. To prove that these guesses are correct we 
need to calculate the sum of squares of the (guessed) degrees of unipotent representations 
which is perhaps the most interesting part of this paper. To do this I find explicit 
formulas (for each irreducible representation $E$ of $W$) of the fake degreess $d_E(q)$
(see no.11). Then I show that the (guessed) degree polynomials can be expressed as linear 
combinations of the $d_E(q)$ with constant coefficients of the form plus or minus $1/2^s$. This 
anticipates the notion of family of representations of the Weyl group and the role of the 
nonabelian Fourier transform, see no.11, (which in this case happens to be abelian.) Here the use 
of the technique of symbols is crucial. 

\subhead 11\endsubhead
I obtained the classification of unipotent representations of $G^F$ (assuming that $G$ is simple of 
exceptional type and $q$ is large enough) in late 1977 (for types $\ne E_8$) and in the spring of 1978 (for 
type $E_8$) when I was already at MIT; this has appeared in 1978 (resp. 1979). In both these papers, as well
as in that in no.10, the notion of fake degree of an irreducible representation $E$ of a Weyl group $W$ is 
considered; it is a certain polynomial $d_E(X)$ with coefficients in $\NN$ whose value at $X=1$ is the 
dimension of $E$. In the second of these papers (on $E_8$) it is
observed that, for any simple $G$, both the set of unipotent representations of $G^F$ and the set
$\Irr W$ can be naturally partitioned into families (the families being indexed
by the same set in both cases). Moreover, to each family $\cf$ one can attach a finite group $\G_\cf$
such that the unipotent representations of $G^F$ in $\cf$ 
are in bijection with the set $M(\G_\cf)$ defined below, the irreducible representations of $W$ in the
family corresponding to $\cf$ are indexed by a subset of $M(\G_\cf)$ and such that the dimension of
unipotent representations in $\cf$ are linear combinations of the fake degress $d_E(q)$ with $E$ running
through the family of irreducible representations of $W$ corresponding to $\cf$. The coefficients in this
linear combination are given by the entries $\{(x,\s),(y,\t)\}$ of a certain matrix indexed by 
$M(\G_\cf)\T M(\G_\cf)$ (see below).

We now define the set $M(\G)$ for any finite group $\G$. It consists of pairs $(x,\s)$ where $x$ is an 
element of $\G$ defined up to conjugacy and $\s\in \Irr Z(x)$ ($Z(x)$  is the centralizer
of $x$ in $\G$). Consider the matrix with entries $\{(x,\s),(y,\t)\}$ indexed by 
$M(\G)\T M(\G)$ given by
$$\{(x,\s),(y,\t)\}=\sum_{g\in\G;xgyg\i=gyg\i x}\fra{\tr(g\i x\i g,\t)\tr(gyg\i,\s)}{|Z(x)||Z(y)|}.$$
This matrix has properties very similar to that of a Fourier transform matrix (it is unitary and
involutive). In fact, when $\G=\FF_2^n$ so that $M(\G)=\FF_2^{2n}$, this is actually a Fourier transform 
matrix.

The fact that this nonabelian generalization of Fourier transform appears in representation theory ($\G$ can
be the symmetric group $S_5$ for $G$ of type $E_8$) was one of the most unexpected things in my research.

\subhead 12\endsubhead
Let $\tW$ be a Coxeter group with canonical set of generators $S$ and length function $\ul:\tW@>>>\NN$. (We 
could take for example $\tW=W$ with $\ul(w)=\dim\co_w-\dim\co_1$, notation of no.7.) Following Iwahori 
(1964) and Bourbaki (1968) we recall the definition of the Hecke algebra $H$ associated to $\tW$. It is an 
algebra over $\ca=\ZZ[v,v\i]$ where $v$ is an indeterminate. As an $\ca$-module, it is free with basis 
$\{T_w;w\in\tW\}$. The (associative) multiplication is defined by the rules
$$T_wT_{w'}=T_{ww'}\text{ if }\ul(ww')=\ul(w)+\ul(w'),$$
$$(T_s+1)(T_s-v^2)=0\text{ if }\ul(s)=1.$$
Note that $T_1$ is the unit element of $H$. Moreover, $\tW$ has a natural partial order. (In the case where 
$\tW=W$, we choose a point $B_0\in\cb$ and for $w\in W$ we set, with notation of no.7,
$o_w=\{B\in\cb;(B_0,B)\in\co_w\}$. We have $y\le w$ if and only if $o_y$ is contained in the closure 
$\bo_w$ of $o_w$.)

\subhead 13\endsubhead
We keep the setup of no.12. Motivated by a result of R. Kilmoyer (1969) in which he described a $v$-analogue 
of the reflection representation of $W$ which was a representation of $H$ in which the generators $T_s$ act 
by particularly simple formulas, I tried (in 1977) to find similar $v$-analogues for other irreducible 
representations of $W$. (I succeeded to do that in only a small number of examples.) This has led me to find
(in the fall of 1978) a new basis of $H$ itself (with $\tW$ as in no.12). First I observed that the map 
$v^nT_w\m\ov{v^nT_w}=v^{-n}T_{w\i}\i$ (with $n\in\ZZ,w\in\tW$) defines a ring involution $h\m\bar h$ of $H$.
Next I showed that for any $w\in\tW$ there is a unique element $C'_w\in H$ such that $\ov{C'_w}=C'_w$ 
and $C'_w=v^{-\ul(w)}\sum_{y\le w}P_{y,w}T_y$ where $P_{y,w}\in\ca$ is a polynomial in $v^2$ of degree 
$\le(\ul(w)-\ul(y)-1)/2$ if $y<w$ and $P_{w,w}=1$. The proof of existence showed also that the polynomials 
$P_{y,w}$ can be computed inductively by a (rather complicated) algorithm. 

Clearly, $\{C'_w;w\in\tW\}$ is an $\ca$-basis of $H$. It has the property that there are many two-sided
ideals of $H$ which are spanned as $\ca$-modules by subsets of $\tW$. The subsets of $\tW$ which appear in 
this way are unions of subsets (called two-sided cells) in a certain partition of $\tW$. 

Consider for example the case where $\tW$ is the Weyl group $W$ of $GL_4(K)$ so that $W$ is the group of 
permutations of $\{1,2,3,4\}$ and let $s_1=(12)$, $s_2=(23)$, $s_3=(34)$. In this case, $P_{y,w}=1$ for any 
$y\le w$ except when $w=s_2s_1s_3s_2$ and $y\in\{1,s_2\}$, or when $w=s_1s_3s_2s_3s_1$ and 
$y\in\{1,s_1,s_3,s_1s_3\}$,  in which case $P_{y,w}=1+v^2$. On the other hand, $\bo_w$ (see no.12) is 
nonsingular except when $w=s_2s_1s_3s_2$ (and its singular locus is $o_1\cup o_{s_2}$) or when 
$w=s_1s_3s_2s_3s_1$ (and its singular locus is $o_1\cup o_{s_1}\cup o_{s_3}\cup o_{s_1s_3}$).

After I told D. Kazhdan about the polynomials $P_{y,w}$ and their apparent connection with the singularities
of $\bo_w$, he suggested a cohomological formula for $P_{y,w}$, assuming that $\tW=W$ and that $P_{y',w}=1$ 
for all $y'$ such that $y<y'\le w$. The construction of the polynomials $P_{y,w}$ and the proof of the 
cohomological formula became part of my joint paper with Kazhdan (published in 1979). That paper also 
contains the conjectural equality 
$$L_w=\sum_{y\le w}(-1)^{\ul(w)+\ul(y)}P_{y,w}(1)M_y\tag a$$
in the theory of highest weight reprepresentations of a simple Lie algebra with Weyl group $W$. ($M_w$ are 
certain Verma modules and $L_w$ is the simple quotient of $M_w$.)

\subhead 14\endsubhead
In 1977, M. Goresky and R. MacPherson introduced intersection cohomology for singular complex algebraic 
varieties which (in the projective case)
satisfied Poincar\'e duality. In early 1979, R. Bott suggested to Kazhdan that in the case 
where $\tW=W$, the $P_{y,w}$ in no.13 might have someting to do with intersection cohomology. I have 
separately arrived at the same conclusion, having attended a talk by MacPherson on intersection cohomology 
at Warwick in 1977; moreover, I knew that (as a consequence of the definition) for any $w\in W$, the 
polynomial $\Pi(X)=\sum_{y\le w}X^{\ul(y)}P_{y,w}(X)$ satisfies $X^{\ul(w)}\Pi(X\i)=\Pi(X)$, which looks 
like a manifestation of Poincar\'e duality. The cohomological
formula for $P_{y,w}$ mentioned in no.13 seemed to be related to a formula in the paper of Goresky and
MacPherson but there was a discrepancy. By talking to MacPherson (in early 1979), Kazhdan and I found out 
that the discrepancy was due to a misprint and also that Deligne has defined intersection cohomology in the 
$l$-adic setting in arbitrary characteristic. After writing to Deligne, we received his letter (in the 
spring of 1979) in which he explained his approach to intersection cohomology using derived categories and 
$l$-adic sheaves. Using Deligne's results, Kazhdan and I were able (in the summer of 1979) to establish the
interpretation of the coefficients of $P_{y,w}$ in no.13 as dimensions of stalks of cohomology sheaves of 
the intersection cohomology complex of $\bo_w$ (in the $l$-adic setting). This appeared in our joint paper 
published in 1980.

The results of this paper were a first step in the proof of the conjecture in no.13(a). The remaining steps
in that proof were achieved in 1981 by Beilinson-Bernstein and Brylinski-Kashiwara. A direct proof of 13(a),
which avoids intersection cohomology, has been found in 2013 by Elias and Williamson, building on work of 
Soergel.

\subhead 15\endsubhead
In July 1979 I gave a talk at AMS conference in Santa Cruz where I stated several conjectures. I will list
here two of them.

(a) I stated a multiplicity formula for the unipotent representations of $G^F$ in the virtual 
representations $\sum_i(-1)^iH^i_c(X_w,\bbq)$ (notation of no.7) in terms of the nonabelian Fourier
transform in no.11. (See no.16.)

(b) I stated a conjecture relating the decomposition of $\Irr W$ into
families (see no.11) with the two-sided cells (see no.13) of $W$. (This was proved by D. Barbasch and 
D. Vogan in 1982 and 1983, based on the solution of the conjecture 13(a).)

\subhead 16\endsubhead
In early 1980 I found a proof of the multiplicity formula in 15(a) in the case where $G$ is of exceptional
type and $q$ is sufficiently large. (This has appeared in a paper in 1980.) In early 1981, while visiting 
the Australian National University, Canberra, I found a proof of the multiplicity formula in 15(a) in the 
case where $G$ is a symplectic or odd orthogonal group and $q$ is sufficiently large; later that year I 
proved the analogous result for even orthogonal groups. (These appeared in papers in 1981 and 1982.) The
methods I used in the 1980 paper were not strong enough to deal with the classical groups treated in the
1981 and 1982 papers. In the last two papers I used the following new approach: instead of studying the 
decomposition into irreducible $G^F$-modules of the $l$-adic cohomology with compact support of $X_w$ in 
no.7 (a very difficult question since $X_w$ is in general not proper) we can try to study the decomposition 
into irreducible $G^F$-modules of the $l$-adic intersection cohomology of the closure $\bX_w$ of $X_w$ in 
$\cb$ (which is better behaved due to the fact that Deligne's theory of weights is applicable so that we can
obtain information about individual cohomology spaces rather than alternating sums of them, which makes a 
crucial difference). The passage from the first problem to the second problem could be controlled at the 
level of Euler characteristic, using the knowledge of the local intersection cohomology of $\bX_w$ (since 
$\bX_w$ looks locally like $\bo_w$, see no.12, it has the same local intersection cohomology as $\bo_w$ 
which, by no.14, is expressible in terms of the polynomials $P_{y,w}$ in no.13). These two papers also rely 
on a method which was new at the time, namely the determination in many cases of the leading coefficients of
the character values of the Hecke algebra.

\subhead 17\endsubhead
From late 1981 to the middle of 1982 I worked on what was to become my 1984 book on the representation
theory of $G^F$ where $G$ is assumed to have connected centre. In this book I obtained the classification of
not necessarily unipotent representations and the computation of their multiplicities in the various 
$\sum_i(-1)^iH^i_c(\tX_w,\bbq)_\th$ (see no.7). For classical groups the classification was already 
essentially known, see no.10, but for exceptional groups it was new for the non-unipotent case and also for 
the unipotent case with $q$ small, see no.11. The multiplicity formulas in the non-unipotent cases and in 
the unipotent case with $q$ small were also new. To study not necessarily unipotent representations, I had 
to study the local and global intersection cohomology of $\bX_w$ with coefficients in certain local systems 
and for this I first had to generalize the results in no.14 to include ``monodromic systems'' on varieties 
like $\bo_w$, see no.12. But the use of intersection cohomology of $\bX_w$ was along the same lines as in 
no.16. Another new ingredient which was not used in the papers in no.16 was the use of the Barbasch-Vogan 
results 15(b); moreover, a result of A. Joseph on Goldie rank representations of Weyl groups was used in the 
proof. (Note that the use of these two ingredients requires the validity of 13(a).)

One of the main results of the book was that any fibre of the map $\Ph$ in no.7 is in a bijection (which 
could be called ``quasi-induction'') with the set of unipotent representations of a smaller $G$ and that the 
multiplicity formulas for the representations in the fibre look similar to those for the corresponding 
unipotent representations. (At the present time there is no a priori proof of this fact.)

The book also contains results for $G$ with not necessarily connected centre. But the case where
$G=Spin_{4k}(\FF_q)$ with $q$ odd could not be treated in the book. It required extensive additional 
computations which I carried out after the book was written, in the summer of 1983, so that in my ICM talk
(1983) I did not have to make any assumption on $G$. (My paper with the details of these computations 
appeared only in 2008.) 

\subhead 18\endsubhead
A unipotent representation $\r$ of $G^F$ is cuspidal (or discrete series) if the following condition is 
satisfied: whenever $w\in W$ is such that $\r$ appears in $\sum_i(-1)^iH^i_c(X_w,\bbq)$, $w$ must be 
elliptic (that is, its conjugacy class does not meet any proper subgroup of $W$ generated by elements of 
length $1$). In some sense, the classification of arbitrary irreducible representations of $G^F$ can be 
reduced to the classification of unipotent cuspidal representations of $G^F$ and of smaller groups. (All 
other irreducible representations are obtained either by quasi-induction, see no.17, at least when $G$ has 
connected centre, or by decomposing some induced representations governed by explicitly known Hecke 
algebras.) In this sense the unipotent cuspidal representation are the most basic of all irreducible 
representations.

Let $\fU_q^0$ be the set of unipotent cuspidal representations of $G^F$ (up to isomorphism). I will describe
a parametrization of the set $\fU_q^0$ which follows from results in my 1984 book and results in one of my 
papers which appeared in 2002; I shall assume that $G$ is simple and split over $\FF_q$ (see no.9).

Let $\r\in\fU_q^0$. There is a unique a pair $(C,\mu)$ where $C$ is a conjugacy class in $W$ and 
$\mu\in\bbq^*$ are such that the three conditions below are satisfied (we denote by $w$ any element of 
minimal length of $C$):

the multiplicity of $\r$ in $H^i_c(X_w,\bbq)$ is $1$ if $i=\ul(w)$ and is $0$ if $i\ne\ul(w)$;

if $z\in W-C$ satisfies $\ul(z)\le\ul(w)$ then the multiplicity of $\r$ in $H^i_c(X_z,\bbq)$ is $0$ for any
$i\in\ZZ$;

the Frobenius map acts on the $\r$-isotypic component of $H^{\ul(w)}_c(X_w,\bbq)$ as multiplication by $\mu$.
\nl
Note that $\mu$ is necessarily a root of $1$ times $q^{\ul(w)/2}$ with $w$ as above. Let $\fS_q$ be the set 
of all pairs $(C,\mu)$ that are attached to some $\r\in\fU_q^0$ as above. The map $\fU_q^0@>>>\fS_q$, 
$\r\m(C,\mu)$, is bijective. (Moreover, two elements of $\fU_q^0$ which have the same associated $\mu$ 
coincide.) 

We now describe the set $\fS_q$ in each case; we will specify $G$ by its type $A_n$, $B_n,\do$, $E_8$. (For
exceptional types we specify an elliptic conjugacy class $C$ in $W$ by the characteristic polynomial $|C|$ of
one of its elements in the reflection representation of $W$, written as a product of cyclotomic polynomials 
$\Ph_d$; we denote by $\th,\sqrt{-1},\z$ a primitive root of $1$ of order $3,4,5$ in $\bbq$.)

$G$ of type $A_n$ $(n\ge1)$: we have $\fS_q=\emp$.

$G$ of type $B_n$ or $C_n$ ($n\ge2$): if $n=k^2+k$ for some integer $k\ge1$, then 
$$\fS_q=\{(C,(-1)^{n/2}q^{k(k+1)(2k+1)/3})\}$$
where $C$ consists of the elements of $W$ which, as permutations of $\{1,2,\do,2n\}$ are a product of cycles
of length $4,8,12,\do,4k$; if $n$ is not of this form, then $\fS_q=\emp$.

$G$ of type $D_n$ ($n\ge4$): if $n=k^2$ for some even integer $k\ge2$, then 
$$\fS_q=\{(C,(-1)^{n/4}q^{2k(k^2-1)/3})\}$$ 
where $C$ consists of the elements of $W$ which, when viewed as permutations of $\{1,2,\do,2n\}$, are a 
product of cycles of length $2,6,10,\do,4k-2$; if $n$ is not of this form, then $\fS_q=\emp$.

$G$ of type $E_6$: $\fS_q$ consists of $(C,\th q^3),(C,\th^2 q^3)$ where $|C|=\Ph_{12}\Ph_3$.

$G$ of type $E_7$: $\fS_q$ consists of $(C,\sqrt{-1}q^{7/2}),(C,-\sqrt{-1}q^{7/2})$ where 
$|C|=\Ph_{18}\Ph_2$.

$G$ of type $E_8$: $\fS_q$ consists of 
$$(C_{30},-\th q^4),(C_{30},-\th^2 q^4),(C_{30},\z q^4),(C_{30},\z^2q^4),(C_{30},\z^3q^4),(C_{30},\z^4q^4),
$$ 
$$(C_{24},\sqrt{-1}q^5),(C_{24},-\sqrt{-1}q^5),(C_{18},\th q^7),(C_{18},\th^2q^7),$$
$$(C_{12},q^{10}), (C'_{12},-q^{11}),(C_6,q^{20}),$$
where 
$$|C_{30}|=\Ph_{30},|C_{24}|=\Ph_{24},|C_{18}|=\Ph_{18}\Ph_6,|C_{12}|=\Ph_{12}^2,|C'_{12}|=\Ph_{12}\Ph_6^2,
|C_6|=\Ph_6^4.$$

$G$ of type $F_4$: $\fS_q$ consists of 
$$(C_{12},\sqrt{-1}q^2),(C_{12},-\sqrt{-1}q^2),(C_{12},\th q^2),(C_{12},\th^2q^2),$$ 
$$(C_8,-q^3),(C_6,q^4),(C_4,q^6),$$
where $|C_{12}|=\Ph_{12}$, $|C_8|=\Ph_8$, $|C_6|=\Ph_6^2$, $|C_4|=\Ph_4^2$.

$G$ of type $G_2$: $\fS_q$ consists of $(C_6,\th q),(C_6,\th^2q),(C_6,-q),(C_3,q^2)$ where $|C_6|=\Ph_6$, 
$|C_3|=\Ph_3$.

\subhead 19\endsubhead
Since the group $G^F$ has a counterpart over $\CC$ one can ask whether the set of unipotent (or unipotent 
cuspidal) representations of $G^F$ has a counterpart over $\CC$. The answer is yes, as we will explain below.

In a paper written in 1977 (which appeared in 1979) I tried to see what happens if, in the definition of the
variety $X_w$ in no.7, one replaces the Frobenius map by conjugation by $t$, a regular semisimple element of
$G$. (Thus, we can define ${}^tY_w=\{B\in\cb;(B,tBt\i)\in\co_w\}$.) These varieties again form a partition 
of $\cb$ into pieces which can be shown to be nonempty and smooth. Now there is no longer an action of an 
interesting group on ${}^tY_w$ but, if $t$ is allowed to vary, the cohomologies of ${}^tY_w$ form 
interesting local systems on the variety of regular semisimple elements in $G$. In this paper I showed that 
from these local systems one can recover information about the unipotent representations $\r$ of $G^F$ (with
$G$ split over $\FF_q$) such that $\r$ appears in the space of functions on the finite set $X_1$. In 
particular, I could show that for such $\r$ the multiplicity of $\r$ in $\sum_i(-1)^iH^i_x(X_w,\bbq)$ is 
independent of $q$ and in fact is expressible in terms of geometry over $\CC$ (since ${}^tY_w$, unlike 
$X_w$, made sense over $\CC$). This paper was for me the beginning of a geometric theory of characters of 
$G$.

\subhead 20\endsubhead
In 1980, I tried to compute in the case where $G=SL_2(K)$, the intersection cohomology complex $\ck$ of $G$ 
with coefficients in the nontrivial local system of rank $1$ on the set of regular semisimple elements given
by the natural double covering of this set. To do this, I noticed that this double covering, when extended 
to the whole of $G$ as the Grothendieck-Springer resolution $\{(g,B)\in G\T\cb;g\in B\}@>>>G$, is a small 
map (in the 
sense of Goresky-MacPherson) which allowed me to perform calculations. These calculations showed that the 
alternating sum of traces of the Frobenius map on the stalks of the cohomology sheaves of $\ck$ (the 
``characteristic function'' of $\ck$) was giving exactly the values of the character of the $q$-dimensional 
(or Steinberg) irreducible representation of $SL_2(\FF_q)$. This was for me a strong indication that for a 
general $G$, the characters of irreducible representations of $G^F$ are closely related to the 
characteristic functions of certain simple perverse sheaves on $G$ and I thus started to look for those 
perverse sheaves. As a first step, I observed that the Grothendieck-Springer resolution was small for 
general $G$ and, as 
an application, I found a new construction of the Springer representations of the Weyl group $W$ in terms of
intersection cohomology which, unlike Springer's original definition, made sense in arbitrary 
characteristic. This appeared in my paper in 1981. At the time (1982) when my 1984 book was written, I had 
the definition of the required collection of simple perverse sheaves (or character sheaves) in the case 
where $G$ was $GL_n$, but I knew that for other groups the analogous collection was not complete and I 
conjectured that it can be completed.

Later, in 1983, I found two definitions of the collection of character sheaves on a general $G$. One, which
appeared in a paper in 1984, was describing the character sheaves as explicit intersection cohomology
complexes on certain subvarieties of $G$; this involved a generalization of the Springer representations of
$W$ and a generalization of the Grothendieck-Springer resolution which satisfied the definition of a small 
map except
that it was not proper in general (but it could still be used). A second definition was given in a series of 
papers on character sheaves which appeared in 1985-1987. It used the idea in no.19 but applied to not 
necessarily regular semisimple elements (that series of papers also contains a proof of the equivalence of 
the two definitions, except for exceptional groups in small characteristic where the proof was completed only
in a 2012 paper). The second definition has the virtue that it is in many ways similar to the constructions 
in no.7; in particular the notion of unipotent character sheaf (analogous to that of unipotent 
representation of $G^F$) is defined.

\subhead 21\endsubhead
We now give the definition of unipotent character sheaves on $G$. For each $w\in W$ we define
$$\p_w:\bY_w=\{(g,B)\in G\T\cb;(B,gBg\i)\in\bco_w\}@>>>G$$ 
by $(g,B)\m g$; here $\bco_w$ is the closure of $\co_w$ in $\cb\T\cb$. Let $\ck_w$ be the direct image under
$\p_w$ of the intersection cohomology complex of $\bY_w$ with coefficients in $\bbq$. By a theorem of 
Beilinson, Bernstein, Deligne and Gabber (1982), $\ck_w$ is a direct sum of simple perverse sheaves (with 
shifts) on $G$. The simple perverse sheaves which appear in this way (for some $w$) are called the unipotent
character sheaves of $w$. Assume that $G$ is simple. A unipotent character sheaf $A$ on $G$ is said to be 
cuspidal if the following condition is satisfied: the set of $g\in G$ such that $A|\{g\}\ne0$ is a union of 
finitely many conjugacy classes (it is then a single conjugacy class).

We now describe the classification of unipotent cuspidal character sheaves of $G$. Let $\fS_1$ be the set of
pairs $(C,\mu_1)$ obtained by replacing each pair $(C,\mu)\in\fS_q$ (see no.18) by $(C,\mu|_{q=1})$. For 
example, if $G$ is of type $A_n$ $(n\ge1$) we have $\fS_1=\emp$. For $G$ of type $B_n$ or $C_n$ ($n\ge2$), 
if $n=k^2+k$ for some integer $k\ge1$, then $\fS_1=\{(C,(-1)^{n/2})\}$ where $C$ is as in no.18 and 
$\fS_1=\emp$ if $n$ is not of this form.
One can show that the set of unipotent cuspidal character sheaves of $G$ is in natural bijection with 
$\fS_1$. Since $\fS_1$ is in natural bijection with $\fS_q$ in no.18 (in type $E_7$ this bijection depends on
a fixed choice of $\sqrt{q}$) we see that the set of unipotent cuspidal character sheaves of $G$ is in 
natural bijection with the set of unipotent cuspidal representations of $G^F$. (A more satisfactory 
explanation of this fact is given in a preprint I posted in 2014.) In this way, the unipotent cuspidal 
character sheaves of $G$ (which can be defined also over $\CC$) appear as limits as $q$ tends to $1$ of 
unipotent cuspidal representations of $G(\FF_q)$.

\subhead 22\endsubhead
The simple Lie group of type $E_8$ (over $\CC$) has an almost mythical status in mathematics. It is the 
simple Lie group for which the dimension divided by the square of its rank is maximum possible (namely 
$\frac{248}{8^2}=4-\frac{1}{8}$) which makes it the most noncommutative simple Lie group. It is the simple 
Lie group for which the dimension of the smallest nontrivial irreducible representation divided by the 
dimension of the group is maximum possible (namely $\frac{248}{248}=1$). It is the only simple Lie group 
such that the centralizer of one of its elements has a nonsolvable group of connected components. It is the 
only simple Lie group $\GG$ of maximal dimension with the following property: there exist conjugacy classes 
$C_2,C_3,C_5$ of elements of order $2,3,5$ respectively such that $\dim C_2+\dim C_3+\dim C_5=2\dim\GG$. (In
fact, $C_2,C_3,C_5$ are uniquely determined; they have dimension $248-240/2,248-240/3,248-240/5$ 
respectively.) It is the only simple Lie group $\GG$ of maximal dimension such that the following set 
$\GG_*$ is nonempty: $\GG_*$ is 
the set of triples $(x_2,x_3,x_5)\in\GG^3$ such that $x_2^2=x_3^3=x_5^5=x_2x_3x_5=1$ and such that 
$\{g\in\GG;gx_2=x_2g,gx_3=x_3g,gx_5=x_5g\}$ is finite. (The fact that $\GG_*\ne\emp$ in type $E_8$ was first
shown by A. V. Borovik in 1989.) Note that $\GG_*$ can be viewed as the set of homomorphisms $A_5@>>>\GG$
(where $A_5$ is the alternating group in $5$ letters) whose image has a finite centralizer.
Let $\GG$ be a simple Lie group of type $E_8$ over $\CC$. The following question was raised by D. D. Frey 
and J.-P. Serre in 1998 (I have learned about it from R. Griess in 2001.) What is the number of orbits of 
the simultaneous conjugation action of $\GG$ on $\GG_*$? Eventually I found that this number is $1$ (this
appeared in my paper in 2003.) It is interesting that the solution of this problem relies on the results on
the representation theory of the group of type $E_8$ over a finite field $\FF_q$, 
explained earlier in this paper.
Frey and Serre have
 shown (in 1998) that $\GG_*=\{(x_2,x_3,x_5)\in C_2\T C_3\T C_5;x_2x_3x_5=1\}$ where $C_n$ are as 
above. Moreover, from the definitions, the $\GG$-action on $\GG_*$ has finite stabilizers.
It is then enough to show that if $\GG$ is replaced by the corresponding group $E_8(\FF_q)$ over 
$\FF_q$ (with $q$ large and not divisible by $2,3,5$) and $C_2,C_3,C_5$ by the corresponding conjugacy 
classes in $E_8(\FF_q)$ (denoted again by $C_2,C_3,C_5$) then 
$$|\{(x_2,x_3,x_5)\in C_2\T C_3\T C_5;x_2x_3x_5=1\}|\tag a$$ 
is approximately equal to $|E_8(\FF_q)|$. By a general result of W. Burnside, the number (a) is equal to
$$\frac{|C_2||C_3||C_5|}{|E_8(\FF_q)|}\sum_\r\frac{\r(g_2)\r(g_3)\r(g_5)}{\r(1)}\tag b$$
where $\r$ runs over the irreducible characters of $E_8(\FF_q)$ and $g_2,g_3,g_5$ are elements of 
$C_2,C_3,C_5$. It is then enough to show that the sum (b) is approximately equal to $q^{248}$. This can be 
done using the classification of irreducible characters of $E_8(\FF_q)$ and the knowledge of their values at 
$g_2,g_3,g_5,1$. While it is too difficult to compute the sum (b) exactly, it is possible to compute it
approximately and this is enough for the desired result. This gives an example of the use of representation 
theory of reductive groups over $\FF_q$ to obtain information on the structure of the corresponding groups 
over $\CC$.

\subhead 23\endsubhead
Another direction of my research is the representation theory of semisimple groups over $p$-adic fields.
In a 1983 paper I found new examples of square integrable representations of such groups associated with
subregular unipotent elements. This relied on some complicated calculations which had the side benefit that
they suggested a way to refine the Deligne-Langlands conjecture (for representations of affine Hecke
algebras) from being a rough parametrization to being an exact parametrization. They also suggested that
the representations of the affine Hecke algebras were intimately related to the geometry of 
``Springer fibres'' in the dual group, in particular that the weight structure of the representations could
be predicted from the geometry of Springer fibres. In my paper written in the summer of 1984, which appeared
in 1985, I 
formulated the idea that equivariant $K$-theory of Springer fibres should be used to construct 
representations of affine Hecke algebras and also that the parameter $q$ of the affine Hecke algebra should 
be interpreted as the generator of the equivariant $K$-theory of a point with respect to the circle group 
action. This idea was further developed in my two papers with Kazhdan (one, written in the fall of 1984 and
one written in the summer of 1985) where the irreducible representations of affine Hecke algebras with equal
parameters were classified, thus providing a proof of the Deligne-Langlands conjecture in the refined form 
stated in my 1983 paper.

In my 1983 paper I also gave a definition of unipotent representations of a simple split adjoint group over a
$p$-adic field and stated a conjectural parametrization for them extending the Deligne-Langlands conjecture.
This parametrization was established in my paper in 1995 where the unipotent representations of finite
reductive groups and the theory of character sheaves come together in a rather unexpected way; I did this by
using character sheaves mixed with equivariant homology to realize representations of various graded
affine Hecke algebras from which the representations of affine Hecke algebras with unequal parameters
could be reconstructed.

\subhead 24\endsubhead
There are several instances where a bar operator similar to the one mentioned in no.13 (in the construction
of the new basis of a Hecke algebra) has been used to construct new bases after its first use in 1978; I 
will mention some of them.

(a) In a 1981 paper I constructed periodic analogues of the polynomials $P_{y,w}$ for affine Weyl groups. In 
my 1999 paper I extended this to a conjecture (proved by R. Bezrukavnikov and I. Mirkovic in 2012) defining
canonical bases in certain equivariant K-theory groups, related to the representation theory of simple Lie 
algebras in characteristic $>0$.

(b) In a 1983 paper I extended the method in no.13 to include Hecke algebras with unequal parameters.

(c) In another 1983 paper (with D. Vogan) we used a bar operator to generalize the polynomials $P_{y,w}$ in
no.13 to the setting of symmetric spaces.

(d) In my 1990 paper I defined a canonical basis in the plus part of the quantized enveloping algebra of
type $A,D$ or $E$ using a bar operator method similar to that in no.13. (Later, M. Kashiwara gave another 
construction of the canonical basis, using the bar operator from my 1990 paper.)

(e) In my 2012 paper with D. Vogan we defined a canonical splitting of $P_{y,w}$ as a sum of two polynomials
(when $y,w$ are involutions in $W$) using a bar operator.
\enddocument